\documentclass[12pt]{amsart}
\usepackage{mathtools, amssymb, bm}

\usepackage{geometry}
\usepackage{biblatex}
\usepackage{color}
\usepackage{amsthm}
\usepackage{lmodern}

\usepackage{hyperref}

\newtheorem{theorem}{Theorem}

\newtheorem{lemma}[theorem]{Lemma}

\theoremstyle{definition}
\newtheorem{remark}[theorem]{Remark}

\numberwithin{theorem}{section}
\numberwithin{equation}{section}

\newcommand{\Z}{\mathbb Z}
\newcommand{\N}{\mathbb N}

\newcommand{\eps}{\varepsilon}

\newcommand{\sqfull}{S^{{\mathstrut\hspace{0.05em}\blacksquare}}_f}

\newcommand{\sqfreeid}{S_{\text{Id}}^{\mathstrut\hspace{0.05em}\square}}

\newcommand{\sqfree}{S_f^{\mathstrut\hspace{0.05em}\square}}

\begin{document}

\title{An average number of square-free values of polynomials}

\author{Watcharakiete Wongcharoenbhorn}

\author{Yotsanan Meemark}

\address{Department of Mathematics and Computer Science, Faculty of Science, Chulalongkorn University, Bangkok, Thailand 10330}

\email{w.wongcharoenbhorn@gmail.com}
\email{yotsanan.m@chula.ac.th}

\keywords{Square-free integers, Riemann hypothesis.}
\subjclass[2020]{Primary 11N32; Secondary 11D79.}

\maketitle

\begin{abstract}
Let $f(X)$ be a polynomial over $\Z$. When $f(X)=\text{Id}(X)$, the identity polynomial, a well-known result states that the number of positive integers $n\leqslant N$ such that $f(n)$ is square-free equals $\frac{N}{\zeta(2)}+O(N^{\frac{1}{2}})$. It is expected that the error term is $O_\eps(N^{\frac{1}{4}+\eps})$ for arbitrarily small $\eps >0$. Usually, it is more difficult to obtain similar order of error term for a higher degree polynomial $f(X)$. We show that the average of error terms, in a weak sense, over polynomials of arbitrary degree, is much smaller than the expected order $O_\eps(N^{\frac{1}{4}+\eps})$.
\end{abstract}

\section{Introduction}
\label{sec:introduction}
For $f(X)\in\Z[X]$, let us denote by $\sqfree(N)$ the number of positive integers $n\leqslant N$ such that $f(n)$ is square-free. We have known as in Tenenbaum \cite{Tenenbaum} that
\begin{equation}
\label{eq: counting-identity}
\sqfreeid(N)=\dfrac{1}{\zeta(2)}N+O(N^{\frac{1}{2}}),
\end{equation}
where we denote by $\text{Id}(X)=X$, the identity polynomial. The order of the error term is expected to be $O_\eps(N^{\frac{1}{4}+\eps})$ for any $\eps>0$. This problem seems to be very difficult. Indeed, the error term in question of order $N^{\frac{1}{2}-\delta}$ for some small $\delta>0$, implies certain quasi-Riemann hypothesis, in which the supremum of the real part of all zeroes of the Riemann zeta function is strictly less than $1$. This was also mentioned in the work of Heath-Brown \cite{square-free}.
\\

Obtaining a similar result on error term as in Equation \eqref{eq: counting-identity} for a higher degree polynomial $f(X)$ is normally more difficult. For instance, Estermann \cite{square-free-Estermann} proved for $f(X)=X^2+h$ with a nonzero integer $h$ that 
$$\sqfree(N)=\prod_{p}\left(1-\dfrac{\rho_h(p^2)}{p^2}\right) N+O(N^{\frac{2}{3}}\log N),$$
where $\rho_h(m):=\#\{n\in\Z/ m\Z: n^2+h\equiv 0\mod{m}\}$ for all $m\geqslant 2$. We see that the order of the main term is the same; however, the error term here is weaker than that of \eqref{eq: counting-identity}. The best known bound is due to Friedlander and Iwaniec \cite{Friedlander-sqfree-quadratic}, in which they obtained $O_\eps(N^{\frac{3}{5}+\eps})$ for arbitrarily small $\varepsilon$. In 2012, Heath-Brown \cite{square-free} further improved the error term for the special case $h=1$ of the above formula to be $O_\varepsilon(N^{\frac{7}{12}+\varepsilon})$. This is, however, still fallen short of $\frac{1}{2}$. In fact, we do not have asymptotic formulae for polynomials of degree more than $3$ unconditionally. Shparlinski \cite{Shparlinski-average-sqfree}, and recently, Browning and Shparlinski \cite{square-free-random} addressed this question on average. The latter work improved the allowable range for polynomials of degree $k$, for any integer  $k\geqslant 4$.\\

To formulate our main result, we denote by $c_f$ the natural density of the counting function of square-free integers corresponding to a polynomial $f(X)$. As a convention, a \textit{natural density} of a set $A\subset \N$ is denoted by
$$\delta(A):=\lim_{N\to\infty} \dfrac{\#(A\cap [1,N])}{N},$$
provided that the limit exists. For each $f(X)$, we expect $c_f$ to be
\begin{equation}
\label{eq:cf}
c_f=\sum_{d\geqslant 1}\dfrac{\mu(d)}{d^2}\rho_f(d^2),
\end{equation}
where $\rho_f(m):=\#\{n\in\Z/m\Z: f(n)\equiv 0\pmod{m}\}$.  From Equation \eqref{eq: counting-identity}, we have that
$$c_{\text{Id}}=\dfrac{1}{\zeta(2)},$$
which satisfies \eqref{eq:cf}. For positive integers $H$ and $k\geqslant 2$, we write $\textbf{a}=(a_0,a_1,\dots,a_k)\in\Z^{k+1}$ and denote
$$\mathcal{F}_k(H):=\{a_0+a_1X+\cdots+a_kX^k\in\Z[X]: \textbf{a}\in\mathcal{B}_k(H)\},$$
where $\mathcal{B}_k(H):=\{\mathbf{a}\in\Z^{k+1}: \gcd(a_0,a_1,\dots,a_k)=1, |a_i|\leqslant H,\text{ for all } i=0,1,\dots,k\}$. Inspired by \cite{square-free-random}, we proved the following theorem.
\begin{theorem}
\label{thm:random-quad-polyn-sqfree}
Let $k$ be any integer greater than $1$. There exists $B=k+1$ such that for all $\eps>0$, $A\geqslant B$ and $H\in[N^B,N^A]$,
$$\dfrac{1}{\#\mathcal{F}_k(H)}\sum_{f\in\mathcal{F}_k(H)} \left(\sqfree(N)-c_{f}N\right)\ll_{A,\eps,k} N^{\eps}.$$
\end{theorem}

This suggests a large cancellation in the average of error terms, in a weak sense, over polynomials in $\mathcal{F}_k(H)$.\\

Numbers that are relatives of square-free integers are ``square-full" integers. An integer is called \textit{square-full} if in its prime factorization, every prime exponent is at least $2$. Denote by $\sqfull(N)$ the number of square-full integers not exceeding $N$. Note that, each square-full integer can be written as $e^2d^3$ for a square-free integer $d$ uniquely. We see for arbitrary $\eps>0$ with $\textbf{a}\in\Z^3$ that
\begin{align*}
\sum_{f\in\mathcal{F}_2(H)}\sqfull(N) &\ll \sum_{d\ll H^{\frac{1}{3}}N^{\frac{2}{3}}}\#\left\{(e,d,n,\textbf{a}):n\leqslant N, \dfrac{a_0+a_1n+a_2n^2}{d^3}=e^2\right\}
\\ &\ll \sum_{d\ll H^{\frac{1}{3}}N^{\frac{2}{3}}} \left(\sqrt{\dfrac{H}{d^3}}+1\right)H^2N \ll H^{\frac{5}{2}}N+H^{\frac{7}{3}}N^{\frac{5}{3}}.
\end{align*}
Upon taking $H\geqslant N^{\frac{5}{2}+\eps}$, we obtain that
$$\dfrac{1}{\#\mathcal{F}_2(H)}\sum_{f\in\mathcal{F}_k(H)}\sqfull(N)\ll_\eps N^\eps,$$
since $\#\mathcal{F}_2(H)\gg H^3$. This suggests that the counting function of square-full values of a quadratic polynomial is really small.

\section{Preliminaries and lemmas}
Throughout this note, we consider ordering polynomials via naive height similar to the work of Browning and Shparlinski \cite{square-free-random}. We shall go through the arguments in parallel to this work. Here, we use some tools from the geometry of numbers. Let
$$\Lambda=\{u_1\textbf{b}_1+u_2\textbf{b}_2+\cdots+u_s\textbf{b}_s:(u_1,u_2,\dots,u_s)\in\Z^s\}$$
be an $s$-dimensional lattice with $s$ linearly independent vectors $\textbf{b}_1,\textbf{b}_2,\dots,\textbf{b}_s$ of integral entries. Denote by $\Delta$ the discriminant of $\Lambda$. Next, we need the following consequence of Schmidt's result in \cite{Schmidt}.

\begin{lemma}
\label{lem:lattice}
Let $\lambda_1\leqslant \lambda_2\leqslant\dots\leqslant\lambda_s$ be the succesive minima of a full rank lattice $\Lambda\subset\Z^s$. Then
$$\#(\Lambda\cap[-H,H]^s)=\dfrac{(2H+1)^s}{\Delta}+O_s\left(\left(\dfrac{H}{\lambda_1}\right)^{s-1}+1\right).$$
\end{lemma}
\begin{proof}
Similar to the proof of Lemma 2.1 in \cite{square-free-random}, we have the asymptotic formula
$$\#(\Lambda\cap[-H,H]^s)=\dfrac{(2H+1)^s}{\Delta}+O\left(\sum_{0\leqslant j\leqslant s-1} \dfrac{H^j}{\lambda_1\lambda_2\cdots \lambda_j}\right),$$
which yields the result.
\end{proof}

Now, let $U_k(m,H,N)$ be the number of solutions to the congruence
$$a_0+a_1n+\cdots+a_kn^k\equiv 0 \pmod m,$$
in the variables
$$(a_0,a_1,\dots, a_k)\in\mathcal{B}_k(H)\text{ 
and } 1\leqslant n\leqslant N.$$

\begin{lemma}[\cite{square-free-random}, Lemma 3.2]
\label{lem:R_f}
For $Q\geqslant 1$, we have
$$\sum_{\frac{Q}{2}<q\leqslant Q} \mu^2(q)U_k(q^2,H,N)\ll_\eps \left(\dfrac{H^{k+1}N}{Q}+NQ+H^k+H^kNQ^{\frac{2}{k+1}}\right)(NQ)^{\eps}.$$
\end{lemma}

\begin{remark}
\label{rm:Uk}
When $\mu^2(q)$ is replaced by $\mu(q)$, we may gain a power saving (from the first term) if we employ the Riemann hypothesis by using the bound on the M\"{o}bius function in Theorem 14.25(C) in \cite{Titchmarsh}. As in the proof in \cite{square-free-random}, we define for $m,n\in\N$ the lattice
$$\Lambda_{m,n}=\{\textbf{a}\in\Z^{k+1}:a_0+a_1n+\cdots+ a_kn^k\equiv 0 \pmod{m}\}.$$
Then, $\Lambda_{m,n}$ is full rank with discriminant $\Delta_{m,n}=m$. By Lemma \ref{lem:lattice}, we obtain
$$U_k(q^2,H,N)=\sum_{n\leqslant N}\left(\dfrac{(2H+1)^{k+1}}{q^2}+O_k\left(\dfrac{H^k}{s(q^2,n)^k}+1\right)\right)$$
where $s(q^2,n)$ is the smallest successive minima of $\Lambda_{q^2,n}$. Now, we consider the main sum above for which we multiply through by $\mu(q)$, and sum over $q\in(Q/2,Q]$, it is 
$$(2H+1)^{k+1}N\sum_{\frac{Q}{2}<q\leqslant Q} \dfrac{\mu(q)}{q^2}\ll_\eps \dfrac{H^{k+1}N}{Q^{\frac{3}{2}-\eps}}$$
by partial summation and the equivalent form of the RH as in Theorem 14.25(C) (\cite{Titchmarsh}).
\end{remark}

\section{Proof of Theorem \ref{thm:random-quad-polyn-sqfree}}
In this section, we prove Theorem \ref{thm:random-quad-polyn-sqfree}. Fix $A\geqslant 1$. We may restrict to the set $\mathcal{F}^*_k(H)$ of irreducible polynomials in $\mathcal F_k(H)$ since $\#(\mathcal{F}_k(H)\backslash\mathcal{F}_k^*(H))\ll_\eps H^{k+\eps}$ by the result of Kuba \cite{Kuba}. Hence, focusing on $\mathcal{F}^*_k(H)$ does not affect the allowable range in our result.\\

We start from the identity
$$\mu^2(n)=\sum_{d^2\mid n}\mu(d).$$
Then we have
$$\sqfree(N)=\sum_{d}\mu(d)\rho_f(d^2,N),$$
where $\rho_f(m,N):=\#\{n\leqslant N: f(n)\equiv 0\pmod m\}$. Note that the sum is finite for if $d\gg \sqrt{HN^k}$ then $d^2$ is larger than $\max_{n\leqslant N,f\in\mathcal F_k(H)} |f(n)|$, hence the summand vanishes for those $d$. Fix $D,E$ that will be specified later, we split the sum above into
\begin{equation}
\label{eq:Sf_split}
\sqfree(N)=M_f(N)+O\left(R^{(1)}_f(N)+R^{(2)}_f(N)\right),
\end{equation}
where
$$M_f(N):=\sum_{d\leqslant D}\mu(d)\rho_f(d^2,N),$$
$$R_f^{(1)}(N):= \sum_{D<d\leqslant E}\mu^2(d)\rho_f(d^2,N) \text{ and }R_f^{(2)}(N):= \sum_{E<d\ll \sqrt{NH^k}}\mu^2(d)\rho_f(d^2,N)$$

Now, we consider the average of $M_f(N)$. Write $\textbf{a}=(a_k,a_{k-1},\ldots,a_0)\in\Z^{k+1}$ and $\textbf{n}=(n^k,n^{k-1},\ldots,1)$. Recall that $N^B\leqslant H\leqslant N^A$, we shall bound each term that is of order $\left(N^{O(1)}\right)^\eps$ by $H^\eps$ for convenience. We have that
\begin{align*}
\sum_{f\in\mathcal F_k(H)} M_f(N) &= \sum_{d\leqslant D} \mu(d)\sum_{f\in\mathcal F_k(H)}\rho_f(d^2,N)
\\ &= \sum_{d\leqslant D} \mu(d)\sum_{\substack{|a_1|,|a_2|,\ldots , |a_k|\leqslant H \\ n\leqslant N}} \#\{|a_0|\leqslant H: \textbf{a}\cdot \textbf{n}\equiv 0\pmod{d^2}\}
\\ &= \sum_{d\leqslant D} \mu(d)NH^k\left(\dfrac{H}{d^2}+O(1)\right).
\end{align*}

Therefore, we obtain
\begin{equation}
\label{eq:avg_Mf}
\sum_{f\in\mathcal{F}_k(H)} M_f(N)= NH^{k+1}\sum_{d\leqslant D} \dfrac{\mu(d)}{d^2}+O(NH^{k}D).
\end{equation}

Now, let us consider the sum of $c_f$ over all $f\in\mathcal{F}_k(H)$. We see from the definition of $c_f$ as in Equation \eqref{eq:cf} that
\begin{align*}
\sum_{f\in\mathcal{F}_k(H)}c_f &=\sum_{f\in\mathcal{F}_k(H)}\sum_{d\geqslant 1}\dfrac{\mu(d)}{d^2}\rho_f(d^2)
\\ &=\sum_{d\leqslant D}\dfrac{\mu(d)}{d^2}\sum_{f\in \mathcal{F}_k(H)}\rho_f(d^2)+\sum_{D<d\ll \sqrt{N^kH}}\dfrac{\mu(d)}{d^2}\sum_{f\in\mathcal{F}_k(H)}\rho_f(d^2).
\end{align*}
Recall that $\rho_f(m):=\#\{n\in\Z/m\Z: f(n)\equiv 0\pmod m\}$. The sum over all $f\in\mathcal{F}_k(H)$ of $\rho_f(d^2)$ is evaluated by
\begin{align*}
\sum_{f\in\mathcal{F}_k(H)}\rho_f(d^2) &=\sum_{\substack{|a_1|,|a_2|,\ldots,|a_k|\leqslant H\\ n\in\Z/d^2\Z}}\#\{|a_0|\leqslant H: \textbf{a}\cdot\textbf{n}\equiv 0\pmod{d^2}\}
\\ &=d^2H^{k}\left(\dfrac{H}{d^2}+O(1)\right)=H^{k+1}+O(H^kd^2).    
\end{align*}
From Corollary 2.3 in \cite{square-free-random}, we have for any square-free $q$ that $\rho_f(q^2)\ll_\eps q^\eps\gcd(\Delta_f, q)$ where $\Delta_f$ is the discriminant of $f$. Therefore, we obtain,
\begin{align*}
\sum_{f\in\mathcal{F}_k(H)} c_f &= H^{k+1}\sum_{d\leqslant D}\dfrac{\mu(d)}{d^2}+O(H^kD)+O\left(H^{k+1}\sum_{d>D}\dfrac{\gcd(\Delta_f,d)}{d^{2-\eps}}\right).
\end{align*}
From Equation (4.7) in \cite{Shparlinski-average-sqfree} we have that
$$\sum_{f\in\mathcal F_k(H)} \dfrac{\gcd(\Delta_f,d)}{d^{2-\eps}}\ll_\eps \dfrac{H^{\eps}}{D^{1-\eps}}.$$
Thus, we obtain upon noting that $H\geqslant N^B$ and $D\ll \sqrt{N^kH}$, that
\begin{equation}
\label{eq:avg_cf}
\sum_{f\in\mathcal{F}_k(H)} c_f=H^{k+1}\sum_{d\leqslant D}\dfrac{\mu(d)}{d^2}+O_{A,\eps}\left(H^kD+\dfrac{H^{k+1+\eps}}{D}\right).
\end{equation}
Let us now recall what we have done so far. By combining Equations \eqref{eq:avg_Mf} and \eqref{eq:avg_cf}, we have for $k\geqslant 2$ that
\begin{align*}
\sum_{f\in\mathcal F_k(H)} &\left(\sqfree(N)-c_fN\right) 
\\ &\ll_{\eps} H^\eps\left(NH^kD+\dfrac{NH^{k+1}}{D}\right)+\sum_{f\in\mathcal{F}_k(H)}R^{(1)}_f(N)+\sum_{f\in\mathcal{F}_k(H)}R^{(2)}_f(N).
\end{align*}
Then, we have from Lemma 2.5 in \cite{square-free-random} and similar arguments,
\begin{align*}
R^{(2)}_f(N)\ll_\eps \left(N^\frac{1}{2}\left(\dfrac{H^\frac{1}{2}N^\frac{k}{2}}{E}+\dfrac{HN^k}{E^2}\right)\right)H^\eps.
\end{align*}
Turning to the remaining term $R^{(1)}_f(N)$. We obtain, by Lemma \ref{lem:R_f} with splitting the sum into dyadic ranges, that
\begin{align*}
\sum_{f\in\mathcal{F}^*_k(H)} R^{(1)}_f(H)\ll_\eps \left(\dfrac{H^{k+1}N}{D}+NE+H^kE+H^kNE^{\frac{2}{k+1}}\right)H^\eps.
\end{align*}
Since $H^k\geqslant N$, this reduces to
\begin{align*}
\sum_{f\in\mathcal{F}^*_k(H)} R^{(1)}_f(H)\ll_\eps \left(\dfrac{H^{k+1}N}{D}+H^kE+H^kNE^{\frac{2}{k+1}}\right)H^\eps.
\end{align*}

Recall that $\#(\mathcal{F}_k(H)\backslash \mathcal{F}^*_k(H))\ll_\eps H^{k+\eps}$. It follows from \eqref{eq:avg_Mf} and \eqref{eq:avg_cf} with the bounds for $R^{(1)}_f(N)$ and $R^{(2)}_f(N)$ that
\begin{align*}
&\dfrac{1}{H^\eps}\sum_{f\in\mathcal{F}_2(H)} \left(\sqfree(N) -c_fN\right)
\\ &\qquad\ll_\eps \left(NH^kD+\dfrac{NH^{k+1}}{D}\right)+\left(N^\frac{1}{2}\left(\dfrac{H^\frac{1}{2}N^\frac{k}{2}}{E}\right)+\dfrac{HN^k}{E^2}\right)H^{k+1}
\\ &\qquad\qquad\qquad\qquad\qquad\qquad\qquad\qquad +\left(\dfrac{H^{k+1}N}{D}+H^kE+H^kNE^{\frac{2}{k+1}}\right)+H^{k}N.
\end{align*}
The term $H^kN$ came from the number of reducible polynomials with each term of order $N$. Whence, upon noting that $\#\mathcal{F}_k(H)\gg H^{k+1}$ and grouping similar error terms, we obtain
$$\dfrac{1}{\#\mathcal F_k(H)}\sum_{f\in\mathcal F_k(H)}\left(\sqfree(N)-c_fN\right)\ll_\eps \Delta H^\eps,$$
where 
$$\Delta:= \dfrac{ND}{H}+\dfrac{N}{D}+\dfrac{H^\frac{1}{2}N^{\frac{k+1}{2}}}{E}+\dfrac{HN^k}{E^2}+\dfrac{E}{H}+\dfrac{NE^{\frac{2}{k+1}}}{H}.$$
We choose $D=H^\frac{1}{2}$ and $E=H^\frac{1}{2}N^\frac{k+1}{2}$. Therefore, we have proved
\begin{theorem}
Fix $A\geqslant 1$ and an integer $k\geqslant 2$. For $N^{k+1}\leqslant H\leqslant N^{A}$,
\begin{align*}
\dfrac{1}{\#\mathcal F_k(H)}&\sum_{f\in\mathcal F_k(H)}\left(\sqfree(N)-c_fN\right) 
\\ &\qquad \ll_{A,\eps,k} \left(1+\dfrac{N}{H^\frac{1}{2}}+\dfrac{N^{\frac{k+1}{2}}}{H^{\frac{1}{2}}}+\dfrac{N^2}{H^{1-\frac{1}{k+1}}}\right)H^{\eps}.
\end{align*}
\end{theorem}

To deduce Theorem \ref{thm:random-quad-polyn-sqfree} we choose $H\geqslant N^{k+1}$ so that each term on the right-hand side, apart from $H^\eps$, is $\ll 1$. This completes the proof.

\begin{remark}
We see that an average of error terms, in a weak sense, exhibits a large cancellation. It is more accessible and seems possible to study the variance than to study the sum of absolute values, i.e., the average of the form
$$\dfrac{1}{\#\mathcal{F}_k(H)}\sum_{f\in\mathcal{F}_k(H)} \left(\sqfree(N)-c_fN\right)^2.$$
We may start as in Equation \eqref{eq:Sf_split} and write $R_f(N):=\sqfree(N)-M_f(N)$. We see that
$\rho_f(d^2,N)=\left(\dfrac{N}{d^2}+O(1)\right)\rho_f(d^2).$
Now, we may determine for fixed $D$
\begin{align*}
M_f(N) &= \sum_{d\leqslant D} \mu(d)\rho_f(d^2,N)
\\ &=N\sum_{d\leqslant D} \dfrac{\mu(d)\rho_f(d^2)}{d^2}+\sum_{d\leqslant D} O(\rho_f(d^2)).
\end{align*}
We have from Lemma 2.4 in \cite{square-free-random}, uniformly over $f\in\mathcal F_k(H)$ that
$$\displaystyle\sum_{d\leqslant D} \rho_f(d^2)\ll_\eps D^{1+\eps}.$$
Whence, we obtain
\begin{align*}
\sum_{f\in\mathcal{F}_k(H)}\left(\sqfree(N)-c_fN\right)^2\ll_\eps \sum_{f\in\mathcal{F}_k(H)}\left(N\sum_{d>D} \dfrac{\mu(d)\rho_f(d^2)}{d^2}\right)^2+\sum_{f\in\mathcal{F}_k(H)} R_f^2(N)+H^{k+1}D^2N^\eps.
\end{align*}
Then, it is seen that when we focus on the first term above
\begin{align*}
\sum_{f\in\mathcal{F}_k(H)}\left(\sum_{d>D}\dfrac{\mu(d)\rho_f(d^2)}{d^2}\right)^2 &= \sum_{d,d_1>D}\dfrac{\mu(d)\mu(d_1)}{d^2d_1^2}\sum_{f\in\mathcal{F}_k(H)} \rho_f(d^2)\rho_f(d_1^2),
\end{align*}
and this is where the first difficulty arises. For a relatively prime pair $(d,d_1)$ we can obtain $\rho_f(d^2)\rho_f(d_1^2)=\rho_f(d^2d_1^2)$ and 
\begin{align*}
\sum_{\substack{d,d_1>D\\ \gcd(d,d_1)=1}} &\dfrac{\mu(d)\mu(d_1)}{d^2d_1^2}\sum_{f\in\mathcal{F}_k(H)}\rho_f(d^2d_1^2) 
\\ \qquad &= \sum_{d>D}\dfrac{\mu(d)}{d^2}\sum_{\substack{d_1>D\\ \gcd(d,d_1)=1}} \dfrac{\mu(d_1)}{d_1^2}\left(H^{k+1}+O(H^kd^2d_1^2\cdot1_{d,d_1\ll HN^{\frac{k}{2}}}(d,d_1))\right),
\end{align*}
where $1_{d,d_1\ll HN^\frac{k}{2}}(d,d_1)=1$ if both $d$ and $d_1$ are $\ll HN^\frac{k}{2}$, and 0 otherwise. For fixed $n\in\N$, the consequence of the Riemann hypothesis (RH), we are to use, is the sum of the type (uniformly for $n\leqslant N^{O(1)}$)
$$\sum_{\substack{d\leqslant D \\ \gcd(d,n)=1}} \mu(d)=O_\eps(D^{\frac{1}{2}+\eps}),$$
which is essentially the adaptation of Theorem 14.25(C) in \cite{Titchmarsh} with the proof of Lemma 2.1 of Vaughan's work (\cite{Vaughan}). Although it is possible to gain a power saving from assuming the RH, it seems difficult to handle the sum when we consider it on non-relatively prime pairs. The same difficulty arises when we bound the term involved $R_f(N)$, although we might again reduce the error term by similar arguments as in Remark \ref{rm:Uk}. It might be possible to hope for 
$$\dfrac{1}{\#\mathcal{F}_k(H)}\sum_{f\in\mathcal{F}_k(H)} \left(\sqfree(N)-c_fN\right)^2\ll_\eps N^{\varpi+\eps},$$
for some $\varpi<1$ so that $\left|\sqfree(N)-c_fN\right|\ll_\eps N^{\frac{\varpi}{2}+\eps}$ on average, under the Riemann hypothesis from these arguments.
\end{remark}

\printbibliography[
heading=bibintoc,
title={References}
]

\end{document}